\documentclass[11pt]{article}
\usepackage[T1]{fontenc}
\usepackage{amsfonts}
\usepackage{amssymb}
\usepackage{amsmath}
\usepackage{latexsym}
\usepackage{graphics}


\newcommand{\ud}{\textrm{d}}   
\newcommand{\bC}{\mathbb{C}}   
\newcommand{\bN}{\mathbb{N}}   
\newcommand{\bR}{\mathbb{R}}   
\newcommand{\bZ}{\mathbb{Z}}   
\newcommand{\proof}{\noindent\textbf{Proof: }}  
\newcommand{\ex}{\noindent\textbf{Example: }}   
\newcommand{\mysize}{4in}      
\newtheorem{thm}{Theorem}[subsection]

\newtheorem{algo}[thm]{Algorithm}
\numberwithin{equation}{section}
\newcounter{excount}

\DeclareMathOperator{\Imz}{Im}   
\DeclareMathOperator{\Rez}{Re}   

\begin{document}

\title{Difficulties in Complex Multiplication and Exponentiation}

\author{ Joshua C. Sasmor \\ \footnotesize{Mathematics Program, Division of Natural and
Health Sciences,} \\ \footnotesize{Seton Hill University,} \\
\footnotesize{Box 502F Seton Hill Drive, Greensburg, PA
15601-1599}\\ \footnotesize{sasmor@setonhill.edu}}




\maketitle






\section{Introduction} \label{sec:intro}
Complex numbers create some of the most beautiful pictures in
mathematics.  Mandelbrot sets, Julia sets, and many other
computer-generated images have become everyday images; hundreds of
books have been published and songs have been written (my favorite
is \cite{jc}).  However, the algorithms which drive these
computations make a major assumption about the complex plane which
can create some visual havoc.  Robert Corless (\cite{co}) asked if
there is a way to compensate for this error.  Unfortunately, we
show this is not possible.

\section{Complex Multiplication and Exponentiation}
\label{sec:cmult}
The multiplication of two complex numbers is rather
straight-forward: for $z_1,z_2 \in \bC$, where $z_1=a+bi$ and
$z_2=c+di$ with $a,b,c,d \in \bR$ we have
\begin{equation} \label{eqn:defofcmult}
z_1 z_2 = (a+bi)(c+di) = (ac-bd) + (ad+bc)i.
\end{equation}
However, this is not computationally efficient for exponentiation;
the computation of $(2+3i)^3$ is not as quick using this
definition:
\begin{eqnarray*}
(2+3i)^3 = (2+3i)(2+3i)(2+3i) & = & ((4-9)+(6+6)i)(2+3i) \\
 & = & (-5+12i)(2+3i) \\
 & = & (-10-36)+(-15+24)i \\
 & = & -46+9i.
\end{eqnarray*}

If we convert $2+3i$ to polar form $z=re^{i\theta}$, we can
compute this exponential much faster.  The conversion equations
between rectangular coordinates $z=a+bi$ and polar form
$z=re^{i\theta}$ are:
\begin{eqnarray}
a = r \cos(\theta) & \quad & b=r \sin(\theta) \label{eqn:PtoR} \\
r=\sqrt{a^2+b^2} & \quad & \tan(\theta) = \frac{b}{a}
\label{eqn:RtoP}
\end{eqnarray}
So we convert $2+3i$ to polar form: $r=\sqrt{2^2+3^2}=\sqrt{13}
\approx 3.60555127\ldots$ and $\tan(\theta)=1.5 \Rightarrow \theta
\approx 0.9828$ or $56.31^{\circ}$.  Now we can use some simple
properties of exponents to get:
\begin{equation}\label{eqn:propofexps}
\left(re^{i\theta}\right)^n=r^n e^{i n \theta} \textrm{ for } n
\in \bN
\end{equation}
So we have \[(2+3i)^3 = \left(\sqrt{13}e^{0.9828i}\right)^3 =
\left(\sqrt{13}\right)^3 e^{3(0.9828)i} \approx
46.87216658e^{2.94838117 i}\] Now convert this back into
rectangular coordinates:
\begin{eqnarray*}
a & = & 46.87216658 \cos(2.94838117) = -46 \\
b & = & 46.87216658 \sin(2.94838117) =9
\end{eqnarray*}
so we have $(2+3i)^3=-46+9i$. Perhaps this polar form conversion
seems more tedious, but consider how much easier it makes
computations with higher powers, like $(2+3i)^{29}$ (which, in
case you're wondering, equals
\\ $-13,833,225,534,613,558-3,190,610,873,034,597i$).

A closer examination of the complex exponential function reveals
an interesting property: $e^z$ is not a one-to-one function! Using
Equations \ref{eqn:PtoR} and \ref{eqn:RtoP}, we see that the
periodic nature of the trigonometric functions cause there to be
multiple (co-terminal) angles which satisfy the formulae.  For
example, $1+i$ can be written as $\displaystyle{\sqrt{2 \
}e^{\frac{\pi}{4}i}}$ and as $\displaystyle{\sqrt{2 \
}e^{\frac{9\pi}{4}i}}$. The angles (referred to as the
\emph{argument} of the complex number; $\arg(z)$) differ by a
multiple of $2\pi$, and any number of the form
$\displaystyle{\sqrt{2 \ }e^{(2k\pi +\frac{\pi}{4})i}}$ with $k\in
\bZ$ is also a polar form of $1+i$.  This means that the inverse
of the function $e^z$ is not a function, but we can work around
this.

The standard definition of the logarithm of $z$ requires that we
specify an interval of length $2\pi$ for the angles.  This is
called a \emph{branch} of the logarithm.  The generalized
logarithmic function is defined by
\begin{equation} \label{eqn:logdefn}
\log(z) = \log(re^{i\theta})=\ln(r) + i \theta
\end{equation}
Now we have a point where this function is undefined: when $z=0$
we would need to compute $\ln(0)$ and to find $\theta$, which is
not well defined for 0. Any subset of the complex plane on which
we can define a (one-to-one) function which acts as an inverse of
$e^z$ is called a branch of the logarithm (see any introductory
text of complex analysis for a more detailed exposition; like
\cite{c}, pages 38-40).  This requires that there is a part of the
plane ``removed'' from the domain of $\log(z)$; this is called the
branch cut.  It is a curve from $z=0$ extending out to infinity;
most often it is a straight line. When we us the principal branch
of the logarithm, the cut is the non-positive real axis. This
means that the angles used for $\theta$ are in the interval
$(-\pi,\pi]$. This is the default branch in almost every computer
algebra system in the world; and is called the \emph{principal
branch} of the logarithm. This means that $\log(z)$ is not
analytic on the whole plane.

\begin{figure}[htp]
\centering \resizebox*{\mysize}{!}{\includegraphics{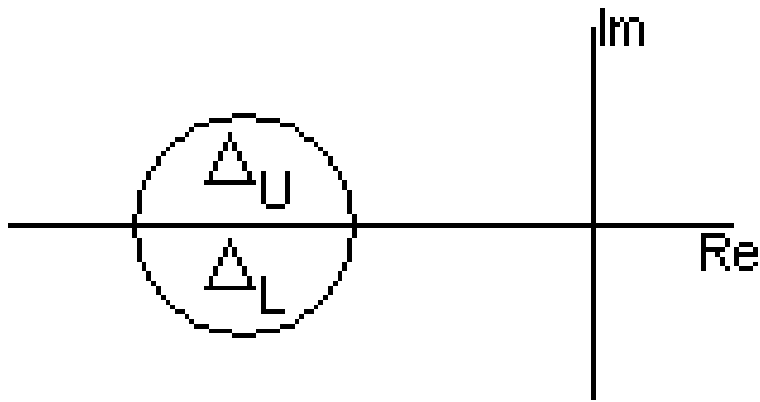}}
\caption[A disk on the negative real axis]{A disk on the negative
real axis} \label{fig:axisdisk}
\end{figure}

When we define $z^n = e^{n \log(z)}$ for $n \in \bN$ we have a
problem: how can the composition of nonanalytic functions become
analytic? It has to do with the behavior of the branch cut. We
will follow the computation through the following chain of
compositions:
\begin{equation} \label{eqn:compchain}
z \mapsto \log(z) \mapsto n \log(z) \mapsto e^{n \log(z)}=z^n
\end{equation}
When we take a small disk $\Delta = \{z:|z-r|<\varepsilon\}$ about
a point $r<0$ in the negative half of the real axis, it is divided
into two halves by the axis (see Figure \ref{fig:axisdisk} on page
\pageref{fig:axisdisk}):
\begin{eqnarray*}
\textrm{the upper half-disk }\Delta_U & = & \Delta \cap \{z|\Imz(z) \geq 0\} \\
\textrm{the lower half-disk }\Delta_L & = & \Delta \cap
\{z|\Imz(z) < 0\}
\end{eqnarray*}

\begin{figure}[htp]
\centering \resizebox*{\mysize}{!}{\includegraphics{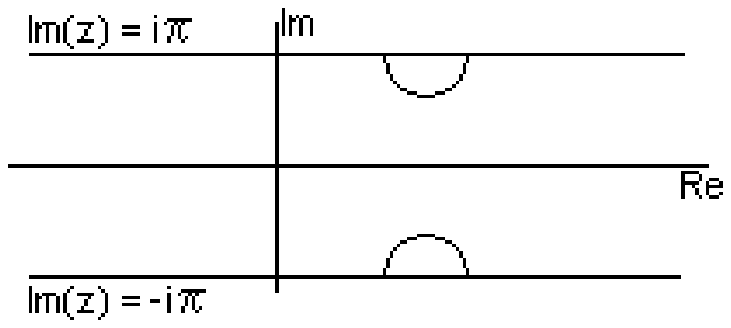}}
\caption[Logarithmic image of the disk]{The image of the disk in
Figure \ref{fig:axisdisk} under the map $\log(z)$}
\label{fig:logdisk}
\end{figure}

These half-disks are transformed by the map $\log(z)$ into regions
about the lines $\Imz(z)=\pi$ and $\Imz(z)=-\pi$.  The images are
almost semi-circular regions, with center at $\log(r)$ and radius
of $\varepsilon$ (see Figure \ref{fig:logdisk} on page
\pageref{fig:logdisk}).

\begin{figure}[htp]
\centering \resizebox*{\mysize}{!}{\includegraphics{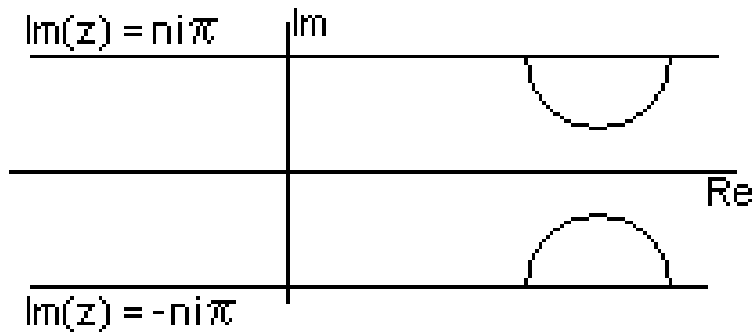}}
\caption[Multiplied logarithmic image of the disk]{The image of
the regions in Figure \ref{fig:logdisk} multiplied by $n$}
\label{fig:nlogdisk}
\end{figure}

Now we have illustrated the first step in the chain in
\ref{eqn:compchain}. Our next step is to multiply these points by
$n \in \bN$.  This expands the entire picture; the horizontal
edges of the half disks are moved to $\Imz(z)=\pm n\pi$, and the
radius will increase by a factor of $n$ as well.  So we have the
image as seen in Figure \ref{fig:nlogdisk} on page
\pageref{fig:nlogdisk}.  Our final step is to exponentiate; the
result will depend on the parity of $n$.  If $n$ is odd, the
result is a disk of radius $\varepsilon^n$ centered around the
point on the negative real axis with magnitude $r^n$ (i.e., the
center is $-r^n$); if $n$ is even, the result is a disk of radius
$\varepsilon^n$ centered around the point $r^n$ (the point on the
positive real axis with magnitude $r^n$).

\begin{figure}[htp]
\centering
\resizebox*{\mysize}{!}{\includegraphics{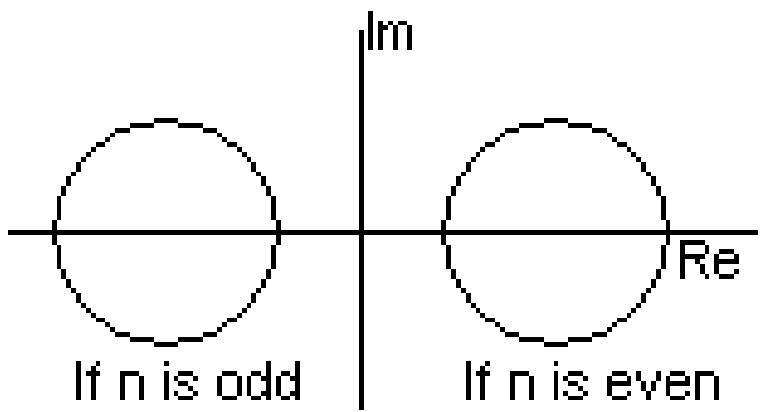}}
\caption[Power map of the disk]{The image of the regions in Figure
\ref{fig:nlogdisk} under the map $e^z$} \label{fig:evenodddisk}
\end{figure}

What is critical in this sequence of maps is the behavior of the
half disks in Figure \ref{fig:nlogdisk} when acted on by the
exponential. When $n$ is an integer, these half disks are aligned
after the exponential map: the arguments of the horizontal
segments of the boundary are $n\pi$ and $-n\pi$ which are mapped
to points on the real axis (see Figure \ref{fig:evenodddisk} on
page \pageref{fig:evenodddisk}). When any other branch cut is
used, the same process can be applied.  The result is that the
function defined by
\begin{equation}
z^n=\left\{\begin{array}{c l} e^{n \log(z)} & \ z \neq 0 \\
0 & \ z=0
\end{array} \right. \end{equation}
is analytic and well defined for $n\in \bN$.

A technical aside: we can consider the infinite family of disks,
with centers on the lines $\Imz(z) = (2k+1)\pi, \  k\in \bZ,$ as
the logarithmic image of the single disk in Figure
\ref{fig:axisdisk}. If we follow this line of reasoning, the
expansion by a factor of $n$ is the critical step; all of the
disks now have centers on lines with $\Imz(z) = nk\pi, \ k \in
\bZ$, and the exponential will map this family of disks (in an
infinite-to-one fashion) onto a single disk as seen in Figure
\ref{fig:evenodddisk} on page \pageref{fig:evenodddisk}.

\section{Non-integer Exponents} \label{sbs:nonint}

What happens when we try to define $z^{\alpha}$ if $\alpha$ is not
an integer? We have the same sequence of composition:
\begin{equation} \label{eqn:compchainalpha}
z \mapsto \log(z) \mapsto \alpha \log(z) \mapsto e^{\alpha
\log(z)}=z^{\alpha}
\end{equation}
but what breaks down?

\begin{figure}[htp]
\centering
\resizebox*{\mysize}{!}{\includegraphics{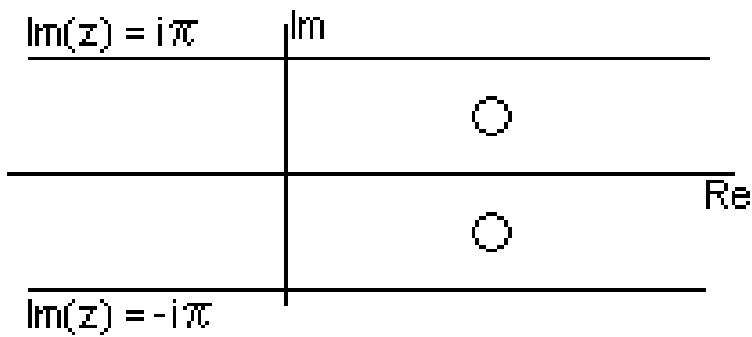}}
\caption[Compressed logarithmic image of a disk]{The images of the
disk in Figure \ref{fig:axisdisk} under the map $\frac{1}{2}
\log(z)$} \label{fig:compresseddisks}
\end{figure}

Let us take a simple example: $z^{\frac{1}{n}}$ for $n\in \bN$.
The process we used for the integer exponents can be paralleled,
but we \emph{must} consider the infinite family of disks mentioned
at the end of Section \ref{sec:cmult}.  When we multiply by
$\alpha = \frac{1}{n}$ in the second step, we will compress $n$ of
these disks to lie between the lines $\Imz(z)=-\pi$ and
$\Imz(z)-\pi$. An example (with $\alpha=2$) is shown in Figure
\ref{fig:compresseddisks} on page \pageref{fig:compresseddisks}.
This results in two distinct square root, and $n$ distinct
$n^{\textrm{th}}$ roots in general, for any non-zero complex
number.  Note that if $n$ is even, then none of the family of
disks will map onto the lines $\Imz(z)=\pm \pi$, and there will be
no problem when we exponentiate.  However if $n$ is odd, then
there are images of the branch cut (the lines $\Imz(z)=\pm n\pi$)
which are compressed onto the lines $\Imz(z)=\pm \pi$. This causes
members of the family of disks to lie on the branch cut, exactly
as in Figure \ref{fig:logdisk}.  Thus when we exponentiate, these
half-disks realign, and the resulting function is analytic (away
from $z=0$).

\begin{figure}[htp]
\centering
\resizebox*{\mysize}{!}{\includegraphics{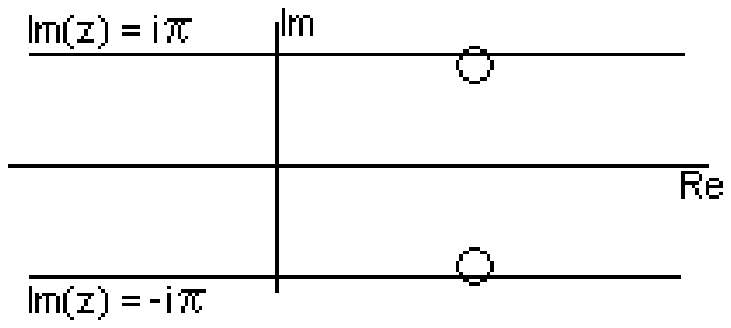}}
\caption[Non-integral compression of a logarithmic image of a
disk]{The images of the disk in Figure \ref{fig:axisdisk} under
the map $\alpha \log(z)$} \label{fig:nonintdisk}
\end{figure}

However, the value $\alpha =\frac{1}{n}$ is a very special case.
What about the more general case? The first step here works
exactly as in Figure \ref{fig:logdisk} on page
\pageref{fig:logdisk}; a disk about the negative real axis is
split into two half-disks and mapped to two half-disks about the
lines $\Imz(z)=\pm \pi$. These regions are expanded by a factor of
$\alpha$, so the lines in figure \ref{fig:nlogdisk} on page
\pageref{fig:nlogdisk} would be labelled $\alpha \pi i$ and
$-\alpha \pi i$.  When we restrict ourselves to the images
contained between the lines $\Imz(z) = \pm \pi$, we have something
similar to that seen in figure \ref{fig:nonintdisk} on page
\pageref{fig:nonintdisk}. Notice that the disks that intersect the
branch cut lines are not bisected, but are divided into nonequal
portions; this is a critical difference.

\begin{figure}[htp]
\centering \resizebox*{2in}{!}{\includegraphics{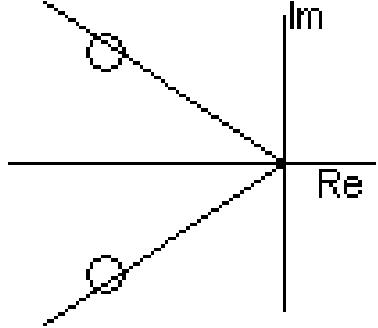}}
\caption[Non-integral power map of the disk]{The image of the disk
from Figure \ref{fig:axisdisk} under the sequence of maps in
Equation \ref{eqn:compchainalpha} on page
\pageref{eqn:compchainalpha}} \label{fig:alphadisk}
\end{figure}

Now we exponentiate, and here is the rub: The images of the rays
\\ $\Imz(z)=\pm \alpha \pi$ do not coincide when $\alpha$ is not an
integer! Depending on the value of $\alpha$, the result is similar
to that seen in figure \ref{fig:alphadisk} on
\pageref{fig:alphadisk}. This causes the function $f(z) =
z^{\alpha}$ not to be analytic on the whole plane; it is not
analytic at the center of our disk, $z=-r$. To be technical again,
we should omit the branch cut from the domain of the function and
say that the function is only analytic on the region $\bC
\setminus (\bR^- \cup \{0\}) = \bC \setminus \{z:\Rez(z) \leq
0,\Imz(z)=0\}$.

However, another definition of an analytic function is a function
which is continuously differentiable on a region $G$, and
$\frac{\ud}{\ud z} (z^{\alpha}) = \alpha z^{\alpha -1}$, so
clearly the function is differentiable.  However, we must be
careful with our choice of arguments (angles) to insure that the
function is \emph{continuously} differentiable.  We will run into
the same argument problem as in the definition of the logarithm,
so we must make a branch cut, leading us to the same domain of
analyticity mentioned in the previous paragraph.

How can we adjust the definition of $z^{\alpha}$ to create a
function that possesses this continuous differentiability?  We
would like to arrange the branch cut so that when we follow the
steps in Equation \ref{eqn:compchainalpha} the branch cuts align
after the exponentiation.  Let us examine a formula for generating
the values of $z^{\alpha}$ and see how to address this problem.

\section{Computation and Formulae} Computationally, what can we
say about $z^{\alpha}$?  How many distinct values in the complex
plane should result?  When $\alpha = \frac{1}{n}$, we have a
direct formula for the $n^{\textrm{th}}$ roots
$z^{\frac{1}{n}}=z^{\alpha}$ (see \cite{c}, for example):
\begin{equation} \label{eqn:nthroots}
z^{\frac{1}{n}}=|z|^{\frac{1}{n}}
e^{\left(\frac{\theta+2k\pi}{n}\right)i} \quad \textrm{where }
\theta = \arg(z), \ 0 \leq k \leq n-1
\end{equation}
which takes on $n$ distinct values in the complex plane.  So if we
were to extend this to non-integer values for $n$, we would have
an equation for a non-integral root of $z$:
\begin{equation} \label{eqn:nonintroots}
z^{\frac{1}{\eta}}=|z|^{\frac{1}{\eta}}
e^{\left(\frac{\theta+2k\pi}{\eta}\right)i} \quad \textrm{where }
\theta = \arg(z)
\end{equation}
But how many values can $k$ range over?  Without loss of
generality, let $\eta = \frac{q}{p}$ with $(p,q)=1$.  Then we can
interpret Equation \ref{eqn:nonintroots} as
\begin{equation} \label{eqn:nonintroots2}
z^{\frac{1}{\eta}}= z^{\frac{p}{q}} = |z|^{\frac{p}{q}}
e^{\left(\frac{(\theta+2k\pi)p}{q}\right)i} \quad \textrm{where }
\theta = \arg(z), \ 0 \leq k \leq q-1
\end{equation}
and we see that this formula will not cause the same value to
appear until we have $q$ distinct values for $z^{\frac{p}{q}}$.
This should make sense, since $z^{\frac{p}{q}} = \left(z^p
\right)^{\frac{1}{q}}$.  When we convert back to the rectangular
coordinates, we see that these $q$ points are exactly the same as
the $q$ points we compute using the $q^{\textrm{th}}$ roots of
$z^p$, but they appear in a different order.

This formula (Equation \ref{eqn:nonintroots2}) can be used to
compute the  $\alpha^{\textrm{th}}$ powers (where $\alpha =
\frac{p}{q}$) of $z$ directly:
\begin{equation} \label{eqn:alphathpower}
z^{\frac{p}{q}} = |z|^{\frac{p}{q}}
e^{\left(\frac{(\theta+2k\pi)p}{q}\right)i} \quad \textrm{where }
\theta = \arg(z), \ 0 \leq k \leq q-1
\end{equation}
However, this is not a one-to-one function; there are $q$ distinct
values. These values have arguments that cover more that an
interval of $2\pi$; in fact these $q$ values have arguments that
range over the interval $[0, 2q\pi)$! This is an important note
and we will return to it after an example.

\ex   Let us compute the $\frac{5}{2}^{\textrm{th}}$ powers of
$1+i$. We will number them $r_0, r_1, \ldots, r_4$.  We have $p=5,
q=2, |1+i|=\sqrt{2}, \theta = \frac{\pi}{4}$. So we can use
Equation \ref{eqn:alphathpower} and we get:
\[r_k=\left(\sqrt{2}\right)^{\frac{2}{5}}
e^{\left(\frac{(\frac{\pi}{4}+2k\pi)2}{5}\right)i}  \textrm{ where
} 0 \leq k \leq 4\] Now we evaluate this to get:
\begin{eqnarray*}
k=0 & \Rightarrow & r_0 = \left(\sqrt{2}\right)^{\frac{2}{5}} e^{\left(\frac{(\frac{\pi}{4})2}{5}\right)i} = \left(\sqrt{2}\right)^{\frac{2}{5}} e^{\left(\frac{\pi}{10}\right)i}\\
k=1 & \Rightarrow & r_1 = \left(\sqrt{2}\right)^{\frac{2}{5}} e^{\left(\frac{(\frac{\pi}{4}+2\pi)2}{5}\right)i} = \left(\sqrt{2}\right)^{\frac{2}{5}} e^{\left(\frac{9\pi}{10}\right)i}\\
k=2 & \Rightarrow & r_2 = \left(\sqrt{2}\right)^{\frac{2}{5}} e^{\left(\frac{(\frac{\pi}{4}+4\pi)2}{5}\right)i} = \left(\sqrt{2}\right)^{\frac{2}{5}} e^{\left(\frac{17\pi}{10}\right)i}\\
k=3 & \Rightarrow & r_3 = \left(\sqrt{2}\right)^{\frac{2}{5}} e^{\left(\frac{(\frac{\pi}{4}+6\pi)2}{5}\right)i} = \left(\sqrt{2}\right)^{\frac{2}{5}} e^{\left(\frac{25\pi}{10}\right)i}\\
k=4 & \Rightarrow & r_4 = \left(\sqrt{2}\right)^{\frac{2}{5}}
e^{\left(\frac{(\frac{\pi}{4}+8\pi)2}{5}\right)i} =
\left(\sqrt{2}\right)^{\frac{2}{5}}
e^{\left(\frac{33\pi}{10}\right)i}
\end{eqnarray*}

Now compare this to the $5^{\textrm{th}}$ roots of $(1+i)^2=2i$,
which we will denote $\rho_0,\ldots,\rho_4$.  We have $|2i|=2,
\theta = \frac{\pi}{2}, n=5$ so by Equation \ref{eqn:nthroots}, we
get:
\[\rho_k = |2|^{\frac{1}{5}}
e^{\left(\frac{\frac{\pi}{2}+2k\pi}{5}\right)i},
\textrm{ where }  0 \leq k \leq 4\] which evaluates to the
following:
\begin{eqnarray*}
k=0 & \Rightarrow & \rho_0 = (2)^{\frac{1}{5}} e^{\left(\frac{\frac{\pi}{2}}{5}\right)i} = (2)^{\frac{1}{5}} e^{\left(\frac{\pi}{10}\right)i}\\
k=1 & \Rightarrow & \rho_1 = (2)^{\frac{1}{5}} e^{\left(\frac{(\frac{\pi}{2}+2\pi)}{5}\right)i} = (2)^{\frac{1}{5}} e^{\left(\frac{\pi}{2}\right)i}\\
k=2 & \Rightarrow & \rho_2 = (2)^{\frac{1}{5}} e^{\left(\frac{(\frac{\pi}{2}+4\pi)}{5}\right)i} = (2)^{\frac{1}{5}} e^{\left(\frac{9\pi}{10}\right)i}\\
k=3 & \Rightarrow & \rho_3 = (2)^{\frac{1}{5}} e^{\left(\frac{(\frac{\pi}{2}+6\pi)}{5}\right)i} = (2)^{\frac{1}{5}} e^{\left(\frac{13\pi}{10}\right)i}\\
k=4 & \Rightarrow & \rho_4 = (2)^{\frac{1}{5}}
e^{\left(\frac{(\frac{\pi}{2}+8\pi)}{5}\right)i} =
(2)^{\frac{1}{5}} e^{\left(\frac{17\pi}{10}\right)i}
\end{eqnarray*}
We can see that there is a one-to-one correspondence between the
values $r_0, \ldots, r_4$ and the values $\rho_0, \ldots, \rho_4$.
\begin{eqnarray*}
\rho_0 = & (2)^{\frac{1}{5}} e^{\left(\frac{\pi}{10}\right)i} & = r_0 \\
\rho_1 = & (2)^{\frac{1}{5}} e^{\left(\frac{\pi}{2}\right)i} & = r_3 \ \left(\textrm{since } \frac{25\pi}{10} \equiv \frac{\pi}{2} \pmod{2\pi} \right) \\
\rho_2 = & (2)^{\frac{1}{5}} e^{\left(\frac{9\pi}{10}\right)i} & = r_1\\
\rho_3 = & (2)^{\frac{1}{5}} e^{\left(\frac{13\pi}{10}\right)i} & = r_4 \ \left(\textrm{since } \frac{33\pi}{10} \equiv \frac{13\pi}{10} \pmod{2\pi} \right)\\
\rho_4 = & (2)^{\frac{1}{5}} e^{\left(\frac{17\pi}{10}\right)i} &
= r_2
\end{eqnarray*}

What we would like to develop is a way to avoid the branch cut
problem seen in Figures \ref{fig:nonintdisk} and
\ref{fig:alphadisk} on page \pageref{fig:alphadisk}.  We would
like to arrange our domain so that the disks which cross the
branch cut always realign (as in Figure \ref{fig:evenodddisk} on
page \pageref{fig:evenodddisk}).  This will require a new domain
for the function $f(z)=z^{\alpha}$.

\section{A New Domain?} \label{sec:newdynspc}

\begin{figure}[htp]
\centering \resizebox*{4in}{!}{\includegraphics{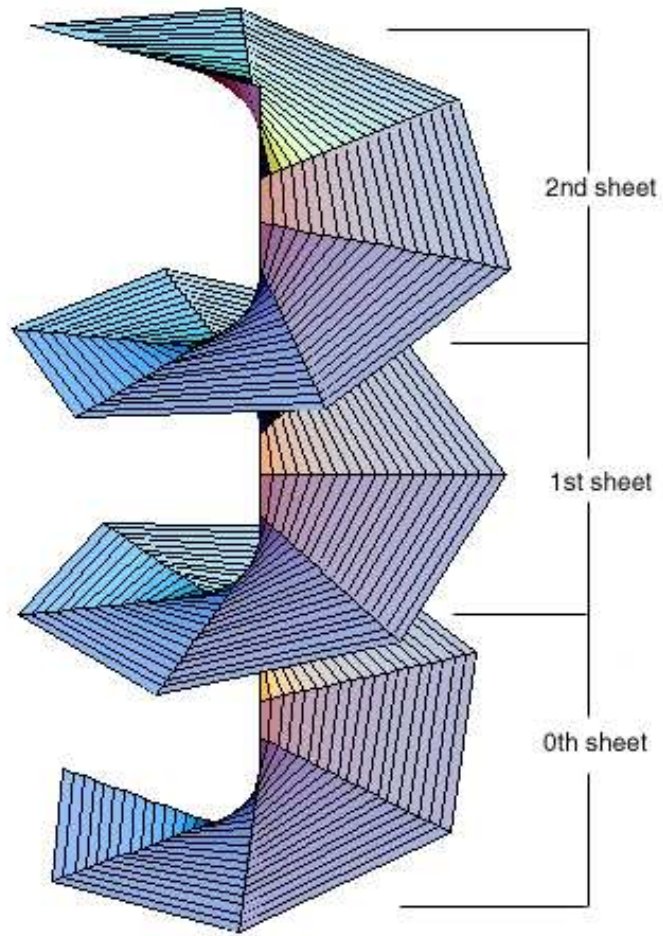}}
\caption[The Riemann surface]{The surface \(W_{\alpha}\)}
\label{fig:spiral}
\end{figure}

For the function $z^{\alpha}+c=e^{(\alpha \log z)}+c$, with
$\alpha = \frac{p}{q}$ a non-integral number, let us redefine the
domain on a more general Riemann surface $W_{\alpha}$. (Note that
this set could also be called $W_q$, but we want to emphasize the
dependence on $\alpha$.) We construct the surface $W_{\alpha}$ by
taking $q$ copies of the plane, all slit along the negative real
axis, each identified by the branch of $\arg (z)$ used on each
sheet. The sheet we label as the $0^{\textrm{th}}$ sheet has
arguments in the interval $(-\pi ,\pi ]$, the sheet we label as
the $+1^{\textrm{st}}$ sheet has arguments in $(\pi ,3\pi]$, and
so on up to the $(q-1)^{\textrm{st}}$ sheet, which has arguments
in the interval $((2q-3)\pi , (2q-1)\pi ]$. In general the
$j^{\textrm{th}}$ sheet has arguments in the interval $((2j-1)\pi
,(2j+1)\pi ]$ for any integer $ 0 \leq j <q$. The sheets are
joined together at the negative real axis so that the argument is
continuously rising as one travels in the positively oriented
direction around the origin, and returning to the
$0^{\textrm{th}}$ sheet after the $(q-1)^{\textrm{st}}$ sheet. In
this fashion, the surface is similar to the Riemann surface for
$z^{1/q}$ with its finite number of sheets. An example is seen in
Figure \ref{fig:spiral} on page \pageref{fig:spiral}.

We must examine the definitions of multiplication and addition in
order to define $F(z) = z^{\alpha} + c$ on this space properly. We
must determine a method to keep track of the sheet we are on when
we operate in this space.

\subsection{Multiplying and Adding in our Domain.} \label{sbs:method}

Define $(z,m)\in (\mathbb{C},\mathbb{Z})$ as $z=re^{2\pi it}=
re^{2\pi i(m+\xi )}$ where $t\in
\left[m-\frac{1}{2},m+\frac{1}{2}\right)$. This assumes that the
branch cut is the negative real axis. We define multiplication
using the same idea as polar coordinates: given
\[z_1 = r_1 e^{2\pi it_1} = r_1 e^{2\pi i(m_1 + \xi _1)}
\qquad \textrm{ and } \qquad z_2 = r_2 e^{2\pi it_2} = r_2 e^{2\pi
i(m_2 + \xi _2)}\] define
\[(z_1,m_1)(z_2,m_2)=(z_1 z_2 , m_*) = (r_1 r_2 e^{2\pi i(t_1 +
t_2)},m_*) \textrm{ where } t_1 + t_2 \in \left[m_* - \frac{1}{2},
m_* + \frac{1}{2}\right)\] Note that $m_*$ may or may not be equal
to $m_1 + m_2$; that depends on the values of $\xi _1$ and $\xi
_2$. However, $m_*$ differs from $m_1 + m_2$ by at most $\pm 1$.
This allows us to define the non-integral powers of $(z,m)$ (i.e.,
$(z,m)^{\alpha} = (z^{\alpha},m_* )$ for the appropriate value of
$m_*$) in terms of the polar coordinates as well. In this context,
keeping track of the sheet information is easy while multiplying.

This means that the Riemann surface $W_{\alpha}$ completely
removes the branch cut discontinuities at 0 and $\infty$
established in the definition of $z^{\alpha}$.  They are replaced
by ramification points of $W_{\alpha}$; i.e., points where
$z^{\alpha}$ is not locally one-to-one.  So we have established a
multiplication that works when $\alpha$ is not an integer.

However, when adding, keeping track of the sheet information is
much harder. We want to define an addition operation
(translation), called $\oplus$ to keep it straight, so that we can
keep track of the sheet information our sum. So suppose that
$(z,m_1) \oplus (c,m_2) = (w,m_*)$, and we define the values of
$w$ and $m_*$ as follows.  Define $w$ as the usual sum of $z+c$
using rectangular coordinate addition.

We use the following cases to determine $m_*$; please note that QI
is the first quadrant of the plane (i.e., the set $\{\Rez(z)>0,
\Imz(z) \geq 0\}$ ), QII is the second quadrant (the set
$\{\Rez(z) \leq 0, \Imz(z) \geq 0\}$), QIII the third quadrant
($\{\Rez(z) \leq 0, \Imz(z)<0\}$), and QIV is the fourth quadrant
($\{\Rez(z)
>0 , \Imz(z)<0\}$). Since we want our operation to agree with the
exponentiation operation, we want the addition to be continuous in
$\theta$, where $z=re^{i \theta}$.  This is referred to as
``Counter-Clockwise Continuity'' or CCC in \cite{coj}.  While this
algorithm seems hard to follow, the idea is simple: if we cross
the negative half of the real axis then we change sheets; if we
cross from QII to QIII, then we go up one sheet, but if we cross
from QIII to QII, then we go down one sheet.

\begin{algo} \label{alg:CCCalg}
To add the point $(c,0)$ to the point $(z,m)$ in $W_{\alpha}$, we
choose the new sheet number $m_*$ as follows:

\begin{enumerate}

\item $\Rez(c)=x_{c}=0$. There are three subcases:

\begin{enumerate}
\item $\Imz(c)=y_{c}=0$. Then $c=0$, $w=z$ and $m_*=m$ and we are
done. \item $\Imz(c)=y_{c}>0$. Then if $z$ is in QI, QII or QIV,
we have \\ $w=z+c$ in the sense of rectangular coordinate addition
and $m_*=m$. However, if $z$ is in QIII, we have two further
subcases:

\begin{enumerate}
\item $|\Imz(c)|\leq |\Imz(z)|$. Then $m_*=m$ \item
$|\Imz(c)|>|\Imz(z)|$. Then $m_*=m-1$. Here is where we actually
cross the branch cut.
\end{enumerate}

\item $\Imz(c)=y_{c}<0$. Then if $z$ is in QI, QIII or QIV, we
have $w=z+c$ in the rectangular coordinate sense and $m_*=m$.
However, if $z$ is in QII, we have two further subcases:

\begin{enumerate}
\item $|\Imz(c)|<|\Imz(z)|$.  Then $m_*=m$. \item $|\Imz(c)| \geq
|\Imz(z)|$. Then $m_*=m+1$. Here we cross the branch cut going the
other way from the crossing above.
\end{enumerate}

\end{enumerate}

\item $\Rez(c)=x_{c}>0$.  Again we have three subcases

\begin{enumerate}
\item $\Imz(c)=y_{c}=0$. Then $m_*=m$ for any $z$. \item
$\Imz(c)=y_{c}>0$. If $z$ is in QI, QII, or QIV then $m_*=m$, but
if $z$ is in QIII there are three subcases:

\begin{enumerate}
\item If $|\Imz(c)| \leq |\Imz(z)|$ then $m_*=m$ \item If
$|\Imz(c)|
> |\Imz(z)|$ and $|\Rez(c)| \geq |\Rez(z)|$ then $m_*=m$ \item If
$|\Imz(c)| > |\Imz(z)|$ and $|\Rez(c)| < |\Rez(z)|$ then
$m_*=m-1$. This subcase, where $c$ moves $z$ up the imaginary
axis, but not over far enough on the real axis to avoid the branch
cut is the only one of these three subcases that crosses the
branch cut.
\end{enumerate}

\item $\Imz(c)=y_{c}<0$. If $z$ is in QI, QIII, or QIV, then
$m_*=m$, but if $z$ is in QII there are three subcases:

\begin{enumerate}
\item If $|\Imz(c)|<|\Imz(z)|$, then $m_*=m$ \item If $|\Imz(c)|
\geq |\Imz(z)|$ and $|\Rez(c)| \geq |\Rez(z)|$ then $m_*=m$ \item
If $|\Imz(c)| \geq |\Imz(z)|$ and $|\Rez(c)| < |\Rez(z)|$ then
$m_*=m+1$.
\end{enumerate}

\end{enumerate}

\item $\Rez(c)=x_{c}<0$. Here also, there are three subcases:

\begin{enumerate}
\item $\Imz(c)=y_{c}=0$. Then $m_*=m$ for any z.

\item $\Imz(c)=y_{c}>0$.
\begin{enumerate}
\item If $z$ is in QI or QII, then $m_*=m$. \item If $z$ is in
QIII, there are two possibilities: Either $|\Imz(c)| \leq
|\Imz(z)|$ in which case $m_*=m$, or $|\Imz(c)|>|\Imz(z)|$ in
which case $m_*=m-1$ \item If $z$ is in QIV, there are three
subcases

\begin{enumerate}
\item If $|\Rez(c)| \leq |\Rez(z)|$ then $m_*=m$ \item If
$|\Rez(c)|
> |\Rez(z)|$ and $|\Imz(c)| \leq |\Imz(z)|$ then $m_*=m$ \item If
$|\Rez(c)|>|\Rez(z)|$ and $|\Imz(c)|>|\Imz(z)|$ then $m_*=m-1$
\end{enumerate}

\end{enumerate}

\item $\Imz(c)=y_{c}<0$

\begin{enumerate}
\item If $z$ is in QIII or QIV then $m_*=m$ \item If $z$ is in QII
the there are two possibilities: Either $|\Imz(c)|<|\Imz(z)|$ in
which case $m_*=m$, or $|\Imz(c)|\geq|\Imz(z)|$ then $m_*=m+1$
\item if $z$ is in QI, then there are three subcases:

\begin{enumerate}
\item if $|\Rez(c)| \leq |\Rez(z)|$ then $m_*=m$ \item if
$|\Rez(c)|>|\Rez(z)|$ and $|\Imz(c)|< |\Imz(z)|$ then $m_*=m$
\item if $|\Rez(c)|>|\Rez(z)|$ and $|\Imz(c)|\geq|\Imz(z)|$ then
$m_*=m+1$
\end{enumerate}

\end{enumerate}

\end{enumerate}

\end{enumerate}
\end{algo}

If we want to generalize to the case where we are adding $(z_1 ,
m_1)$ and $(z_2 , m_2)$, the process is nearly identical. We wish
to preserve the idea from addition of vectors that the sum of two
vectors lies ``between'' the vectors (i.e., the parallelogram
law). The major change in the method above is that $m_*$ is now
based on $\lfloor\frac{m_1 + m_2}{2}\rfloor$ and then adjusted by
$\pm 1$ depending on the cases above. The only main concern is
that the point 0 is on every sheet. Hence if $z_{1} + z_{2} = 0$
then the sheet is irrelevant. By default, the sheet should be left
as $\lfloor\frac{m_1 + m_2}{2}\rfloor$, in order to simplify the
definition. This now allows us to add any two complex numbers in
the Riemann surface setting.  We will define neighborhoods of zero
topologically: any ball of radius $\varepsilon$ around zero,
covering all sheets, is a neighborhood of zero.

The major change is that addition is now a noncommutative
operation! It is order dependent, as we see in the following
example.

\addtocounter{excount}{1}

\textbf{Example \theexcount :} Suppose we are working on a
4-sheeted space, like $W_{15/4}$. Let $(z_1 , m_1) = (-2+i,2)$ and
$(z_2 , m_2) = (-1-3i,2)$. Converting these to polar coordinates,
we get $z_1 \approx \sqrt{5} e^{0.852416 \pi i}$ and $z_2 \approx
\sqrt{10} e^{0.602416 \pi i}$. So when we compute $(z_1 + z_2,
m_*)$ we get $z_1 + z_2 = -3-2i$ and $m_* = 3$ (by case 3(c)(ii)
above).  However, when we compute $(z_2 + z_1, m_*)$ we get $z_2 +
z_1 = -3-2i$ and $m_* = 2$ (by case 3(b)(ii) above). $\Box$

We run into a big problem: this operation is not continuous.


\textbf{Counter-example \theexcount :}  Let us consider an
$\varepsilon$ -neighborhood $N$ of the point $z_0 =
(\frac{-i}{2},1)$, under the map $z \mapsto z+i$, on a three
sheeted surface like $W_{\frac{8}{3}}$.  Then when we translate,
using Algorithm \ref{alg:CCCalg}, the points in $N \cap
\{\Rez(z)<0\}$ are shifted to sheet 0, while the points in $N \cap
\{\Rez(z) \geq 0\}$ remain on sheet 1. Hence $N$ is ``sheared''
into two half-disks: one containing the point $(\frac{-i}{2},1)$,
and one containing the point $(\frac{-i}{2},0)$ in its closure.
$\Box$

This algorithm is specifically constructed to match the choice of
the negative real axis as the branch cut.  This is not the only
way to construct the surface $W_{\alpha}$ and still maintain
continuity in the argument, $t$, in our construction.  We could
use a branch cut along the positive real axis, or along any curve
between the points 0 and $\infty$, just as any of these curves
define a valid branch of the logarithm.  But this will not change
the continuity problem.  In fact, any rigid translation on
$W_{\alpha}$ is inherently discontinuous.  Once we decide that our
domain is a ramified Riemann surface like $W_\alpha$, translation
becomes discontinuous at the ramification point (0 in our case).

\begin{thm} \label{thm:noctsadd}
Let $C$ be the branch cut from 0 to $\infty$ of the surface
$W_{\alpha}$. For any $c \in W_{\alpha}$, there exists a point
$z_0$ and a neighborhood $N = N(z_0)$ such that the set $N+c =
\{z+c : z \in N\}$ is not a connected set (using a translation
defined similarly to Algorithm \ref{alg:CCCalg}, with $C$ as the
branch cut). In particular, we may choose $z_0$ such that $0 \in
\{z_0 +tc, \, 0 \leq t \leq 1 \}$.
\end{thm}

\proof{} Fix $\alpha$ and $c$.  Choose $z_0$ such that $0 \in
\{z_0 +tc, \, 0 \leq t \leq 1 \}$. Then the ray from the origin
through $c$ divides $N=N(z_0)$ into two subsets $N_L$ and $N_U$,
each of which is a half-disk.  When $N$ is translated, points in
one of the half-disks change sheets, while points in the other
half-disk do not change sheets.  Regardless of which sheet the
point $z_0$ is mapped onto, half of $N$ will be separated from
$z_0$.  If we approach $z_0$ in $N_L$, then for the sake of
continuity, $z_0$ should be mapped onto the sheet containing
$N_L$; but the same argument holds true for $N_U$. But the point
$z_0$ is itself only mapped onto one sheet. Hence the translation
map is not continuous.  \(\Box\)

This answers a question raised by Robert Corless in his E.C.C.A.D.
presentation \cite{co}: ``Can a Riemann surface variable be coded?
What will the operations be on it?''  (this was also discussed in
detail in the Appendix A to E. Kaltofen's paper \cite{k}.)
Unfortunately, we are answering ``No, a Riemann surface variable
cannot be coded. There is no continuous addition operation on such
a variable.''

We might also choose a different algorithm for choosing $m_*$. For
example, we can define $m_*$ based on the sign of $\Imz(c)$:

\begin{algo} \label{alg:signalgo}
To add the point $(c,m_2)$ to the point $(z,m_1)$ in $W_{\alpha}$,
we choose the new sheet number $m_*$ as follows:
\begin{enumerate}
\item If $\Imz(c)>0$ then $m_* = m_1+m_2+1$. \item If $\Imz(c)=0$
then $m_* = m_1+m_2$. \item If $\Imz(c)<0$ then $m_* = m_1+m_2-1$.
\end{enumerate}
\end{algo}

But this algorithm for choosing sheets is not continuous either.

\addtocounter{excount}{1} \textbf{Counter-example \theexcount :}
Fix $\alpha = \frac{5}{2}$ and $c=(i,0)$, so $W_{\alpha}$ has two
sheets, with branch cut along the negative real axis.  Consider an
$\varepsilon$-neighborhood $N$ of the point $(-1,0)$. Now half of
$N$ is on sheet 0 and half of $N$ is on sheet 1.  Hence when we
add $i$ using Algorithm \ref{alg:signalgo}, the neighborhood $N$
is ``sheared'' into two parts, since the points on sheet 1 are
moved to sheet 0 and vice-versa.  $N$ will become two half-disks
centered about the points $(i,0)$ and $(i,1)$.

By defining the function $z^{\alpha} = e^{\alpha \log (z)}$ we
have two points in $\overline{\mathbb{C}}$ where $F$ is not
conformal, namely $0$ and $\infty$.  This is caused by the
branching of $\log(z)$ and the failure of the branch cuts to line
back up after multiplying by $\alpha$ (recall the Figures
\ref{fig:nonintdisk} and \ref{fig:alphadisk}). Any translation by
$c$ cannot prevent these points from staying ramified.  Hence
ramification at zero implies either $c$ is also a ramification
point or $z \mapsto z+c$ is not continuous. Since ramification at
$c$ requires infinitely many ramification points (since $c$ is
chosen arbitrarily), this isn't a good plan.

\section{Looking Forward}

We come to the following negative result: the operation of
addition is incompatible with the operation of exponentiation for
non-integer exponents. So in order to study the fractals generated
by functions of this form, we must redefine them on the plane and
make sure that the definitions match the definitions for the
integer case as closely as possible. We have examined these sets
in a separate article (\cite{s3}).

\begin{figure}[h!tp]
\centering \resizebox*{\mysize}{!}{\includegraphics{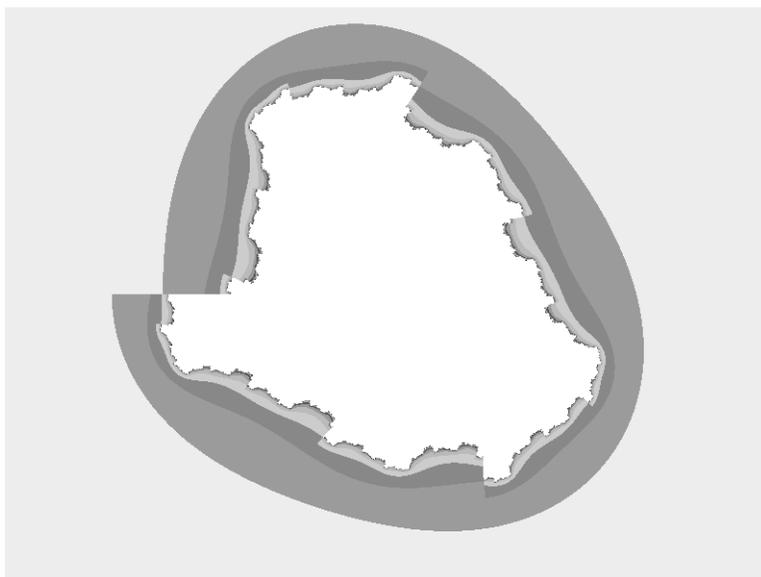}}
\caption[Julia set for $f(z)=z^{2.5}+\frac{1}{2}i$]{The Julia set
for the function $f(z)=z^{2.5}+\frac{1}{2} i$} \label{fig:jset2p5}
\end{figure}

Many of the characterizations of fractals for polynomials will
fail to carry over to these functions and that the discontinuity
we observed in Section \ref{sec:cmult} becomes more important than
ever.  This discontinuity is observed when an algorithm for
generating fractals is given non-integer exponents. Figure
\ref{fig:jset2p5}, on page \pageref{fig:jset2p5}, is an example of
the Julia set for the function $f(z)=z^{2.5}+\frac{1}{2} i$.  This
image was generated using the Fractint program (version 19.6)
found at the web site given in \cite{fr} (using the fractal type
\verb|julzpower|), and the calculations were carried out using
only the principal branch of the logarithm. The parameters for the
\verb|julzpower| fractal type are: the real and imaginary parts of
the parameter ($0+0.5i$); the real and imaginary parts of the
exponent ($2.5 + 0i$); the bailout test and the bailout value,
both of which are set to 0 for the defaults (modulus for the test
and 4 for the value).

\clearpage



\end{document}